\documentclass[10pt]
{article}

\usepackage{amsthm}
\usepackage[margin=30mm]{geometry}
\usepackage{amsmath}%
\usepackage{amsfonts}%
\usepackage{amssymb}%
\usepackage{amscd}
\usepackage{graphicx}
\usepackage{color}
\usepackage{xcolor}
\usepackage{bbm}
\usepackage[all,cmtip]{xy}
\usepackage{rotating}
\usepackage{url}
\usepackage{currfile}
\usepackage{subcaption}
\usepackage{enumerate}
\usepackage{enumitem}
\usepackage{wrapfig}
\usepackage{eurosym}

\definecolor{grau}{rgb}{0.1,0.1,0.1}
\usepackage[colorlinks=true, urlbordercolor=grau, linkcolor=grau, citecolor=grau, urlcolor=grau]{hyperref}
\usepackage{breakurl}

\usepackage{tikz}
\tikzset{node distance=2cm, auto}
\usetikzlibrary{math}

\usepackage{bibspacing}

\newtheorem{theorem}{Theorem}[section]
\newtheorem{proposition}[theorem]{Proposition}
\newtheorem{lemma}[theorem]{Lemma}
\newtheorem{corollary}[theorem]{Corollary}

\newtheorem{question}[theorem]{Question}

\theoremstyle{definition}

\newtheorem{remark}[theorem]{Remark}
\newtheorem{remarks}[theorem]{Remarks}

\newcommand{\vrecica}{Vre\'{c}ica}
\newcommand{\zivaljevic}{\v{Z}ivaljevi\'{c}}

\newcommand{\RR}{\mathbbm{R}} 
\newcommand{\QQ}{\mathbbm{Q}} 
\newcommand{\ZZ}{\mathbbm{Z}} 
\newcommand{\CC}{\mathbbm{C}} 
\DeclareMathOperator*{\wo}{\backslash}

\newcommand{\textdef}[1]{\textnormal{\textit{#1}}}
\newcommand{\one}{\mathbbm{1}} 

\newcommand{\impl}{\Rightarrow} 

\newcommand{\eps}{\varepsilon}

\DeclareMathOperator{\Tr}{Tr}

\newenvironment{psmallmatrix}
  {\left(\begin{smallmatrix}}
  {\end{smallmatrix}\right)}

\newcommand{\mysquare}{%
\begin{tikzpicture}%
\tikzmath{\x = 1.3;\y = \x;}
\draw (0,0) -- (\x ex,0);%
\draw (\x ex,0) -- (\x ex,\y ex);%
\draw (\x ex,\y ex) -- (0,\y ex);%
\draw (0,\y ex) -- (0,0);%
\end{tikzpicture}%
}
\newcommand{\rectangle}{%
\begin{tikzpicture}%
\tikzmath{\x = 1.8;\y = 1.1;}
\draw (0,0) -- (\x ex,0);%
\draw (\x ex,0) -- (\x ex,\y ex);%
\draw (\x ex,\y ex) -- (0,\y ex);%
\draw (0,\y ex) -- (0,0);%
\end{tikzpicture}%
}
\newcommand{\trapezoid}{%
\begin{tikzpicture}%
\tikzmath{\x = 1.8;\y = 1.2;\s=0.5; \xx = \s; \xxx = \x - \s;}
\draw (0,0) -- (\x ex,0);%
\draw (\x ex,0) -- (\xxx ex,\y ex);%
\draw (\xxx ex,\y ex) -- (\xx ex,\y ex);%
\draw (\xx ex,\y ex) -- (0,0);%
\end{tikzpicture}%
}

\newcommand{\C}{\mathcal{J}}
\newcommand{\Cconv}{\C_{\textnormal{conv}}}
\newcommand{\Cconvn}[1]{\Cconv^{#1}}
\newcommand{\CconvZero}{\Cconvn{0}}
\newcommand{\CconvOne}{\Cconvn{1}}

\newcommand{\Cn}[1]{\C^{#1}}
\newcommand{\CZero}{\Cn{0}}
\newcommand{\COne}{\Cn{1}}
\newcommand{\COnePW}{\Cn{1}_{\textnormal{pw}}}
\newcommand{\CInf}{\Cn{\infty}}
\newcommand{\Q}{\mathcal{Q}}
\newcommand{\Qsq}{\Q_{\mysquare}}
\newcommand{\Qrect}{\Q_{\rectangle}}
\newcommand{\Qtrap}{\Q_{\trapezoid}}
\newcommand{\Qcirc}{\Q_{\bigcirc}}

\begin{document}

\title{Quadrilaterals inscribed in convex curves}

\author{%
Benjamin Matschke%
\setcounter{footnote}{-1}%
\\
\small Boston University\\
\small matschke@bu.edu
}
\date{\today}

\maketitle

\begin{abstract}
We classify the set of quadrilaterals that can be inscribed in convex Jordan curves, in the continuous as well as in the smooth case.%
\footnote{The author was kindly notified by Arseniy Akopyan and Sergey Avvakumov about their recent article~\cite{AkopyanAvvakumov17CyclicQuads}, which has considerable overlap with the present one. As they had already submitted their preprint, it was too late to merge them. Certainly, their results have full priority; in particular Theorem~\ref{thmMainCircularQuads} appeared already in their work as well as (as a corollary) the first part of Theorem~\ref{thmMainTrapezoids} in the special case of rectangles.}
This answers a question of Makeev in the special case of convex curves.
The difficulty of this problem comes from the fact that standard topological arguments to prove the existence of solutions do not apply here due to the lack of sufficient symmetry. 
Instead, the proof makes use of an area argument of Karasev and Tao, which we furthermore simplify and elaborate on.
The continuous case requires an additional analysis of the singular points, and a small miracle, which then extends to show that the problems of inscribing isosceles trapezoids in smooth curves and in piecewise $C^1$ curves are equivalent.
\end{abstract}

%
%
%
%

\renewcommand\contentsname{Contents}
\setcounter{tocdepth}{1}

\section{Introduction}

A Jordan curve is a simple closed curve in the plane, i.e.\ 
an injective continuous map $\gamma:S^1\to\RR^2$.
In~1911, Toeplitz~\cite{Toe11aufgabenDerAnalysisSitus} announced to have proved that any convex Jordan curve contains the four vertices of a square -- a so-called inscribed square -- and he asked whether the same property holds for arbitrary Jordan curves.
This became the famous Inscribed Square Problem, also known as the Square Peg Problem or as Toeplitz' Conjecture.
So far it has been answered in the affirmative only in special cases~\cite{Emc13squarePeg1,Hebbert14InscribedSquaresAndKinematicGeometry,
Emc15squarePeg2,Emc16squarePeg3,Zin21konvexeGebilde, Shn44geomPropClosedCurves, Chr50kvadrat,Jer61inscribedSquares, Str89inscribedSquares,
NielsenWright95rectanglesInscribedInSymmetricContinua,
Mak95quadsInscribedInClosedCurve, VrZi08fultonMacPhersonCompCyclohedraAndPolygonalPegProblem, Pak09discreteAndPolyhedralGeometry, SaMa09inscribedSquaresDigitalPlane, SaMa11inscribedSquares, CDM11squarePeg, Mat09squarePeg,Mat11phd,OPT13noteOnToeplitzSquareProblem,vanHeijst14masterThesis,Tao17integrationApproachToToeplitzProblem}.

More generally, we say that a Jordan curve $\gamma$ \textdef{inscribes} a quadrilateral $Q$ if there is an orientation-preserving similarity transformation that sends all four vertices of $Q$ into the image of $\gamma$.
Thus Toeplitz proved that convex Jordan curves inscribe squares.

It is natural to ask whether they inscribe more general quadrilaterals as well. 
This is methodology-wise a highly interesting question for the following reason:
Almost all approaches up to today (with few exceptions, Tao~\cite{Tao17integrationApproachToToeplitzProblem}; and for more general circular quadrilaterals see also Karasev~\cite{KaVo10MakeevsConj}, and for rectangles of aspect ratio~$\sqrt{3}$ see~\cite{Mat11phd}) prove the existence of inscribed squares via more or less directly proving topologically that the number of inscribed squares is \emph{odd} when counted with appropriate multiplicities, and thus never zero.
Any other quadrilateral turns out to be inscribed an even number of times (or zero times when counted with appropriate signs) due to their smaller symmetry group, and thus the topological approach does not extend to quadrilaterals that are not squares.

A \textdef{circular quadrilateral} is a quadrilateral that has a circumcircle.
An \textdef{isosceles trapezoid} is a trapezoid that has a circumcircle.
Let $\Qsq$, $\Qrect$, $\Qtrap$, $\Qcirc$ denote the sets of squares, rectangles, isosceles trapezoids and circular quadrilaterals, respectively.
Clearly, $\Qsq\subset \Qrect\subset \Qtrap\subset \Qcirc$.

If $\C$ is a class of Jordan curves and $\Q$ a set of quadrilaterals, we say that $\C$ inscribes $\Q$ if any curve $\gamma\in\C$ inscribes each quadrilateral $Q\in\Q$.

Let $\Cn{k}$ denote $k$-times continuously differentiable Jordan curves that are regular if $k\geq 1$.
Makeev~\cite{Mak95quadsInscribedInClosedCurve} asked: Does $\CZero$ inscribe $\Qcirc$? One restricts to $\Qcirc$ clearly because the only quadrilaterals that are inscribable in circles are circular.
Quite likely Makeev meant $\COne$ instead of $\CZero$ (compare with Makeev~\cite{Mak05quadsInscribedInCurveAndVerticesOfCurve}), as it turns out that for example the only quadrilaterals that can be inscribed in arbitrarily thin triangles are isosceles trapezoids, as observed by Pak~\cite[Ex. 5.16]{Pak09discreteAndPolyhedralGeometry}.
In any case one arrives at two natural questions.

\begin{question}[Continuous case]
\label{quJ0inscribesQtrap}
Does $\CZero$ inscribe $\Qtrap$?
\end{question}

\begin{question}[Smooth case]
\label{quJ1inscribesQcirc}
Does $\COne$ inscribes $\Qcirc$?
\end{question}

Makeev~\cite{Mak95quadsInscribedInClosedCurve} managed to answer Question~\ref{quJ1inscribesQcirc} in the affirmative in the special case of star-shaped $C^2$-curves that intersect every circle at most $4$ times, see also Makeev~\cite{Mak05quadsInscribedInCurveAndVerticesOfCurve} for a version of that.
To underline the difficulty of both questions, 
the author~\cite{Mat12surveyOnSquarePeg} had put \euro 100 on the weaker problem of whether or not $\CInf$ inscribes $\Qrect$. 

\begin{remark}[Updates] This prize was recently earned by Greene and Lobb~\cite{GreeneLobb20rectangles}, 
and moreover they just announced in~\cite{GreeneLobb20circularQuads} a positive answer for Question~\ref{quJ1inscribesQcirc}.
Their proofs are based on symplectic topology, in particular using the minimum Maslov number of Lagrangian tori in~$\CC^2$ for the general case.
In combination with Theorem~\ref{thmReductionOfCOnePWtoCInf_forInscribingQtrap} this allows us to also provide partial positive answers to Questions~\ref{quJ0inscribesQtrap} and~\ref{quJ1inscribesQcirc} for the class $\COnePW$ of piecewise $C^1$ Jordan curves without cusps; see Corollaries~\ref{corJ1pw_inscribes_trap} and~\ref{corJ1pw_with_large_angles_inscribes_circ}. 
\end{remark}


In the current paper, we answer both questions in the affirmative in the case of convex curves.

\begin{theorem}[Continuous case]
\label{thmMainTrapezoids}
The class $\CconvZero$ of (continuous) convex Jordan curves inscribes the set $\Qtrap$ of isosceles trapezoids.
Moreover, $\Qtrap$ is the largest possible such set of quadrilaterals.
\end{theorem}

\begin{theorem}[Smooth case]
\label{thmMainCircularQuads}
The class $\CconvOne$ of differentiable convex Jordan curves inscribes the set $\Qcirc$ of circular quadrilaterals.
Moreover, $\Qcirc$ is the largest possible such set of quadrilaterals.
\end{theorem}




\paragraph{A common generalization.}

The above two theorems state the inscribability of $\Qcirc$ in $\CconvOne$, and of $\Qtrap$ in $\CconvZero$.
Additionally we know that for each quadrilateral not in $\Qtrap$ there is a curve in $\CZero$ that does not inscribe it.
Nonetheless, we may ask for natural sufficient criteria for when a circular quadrilateral can be inscribed into a (continuous) convex Jordan curve.
One positive answer is given in the following Theorem~\ref{thmMainExtension}.

\begin{figure}[htb]
  \centering
  \begin{minipage}[b]{0.45\textwidth}
	\centering
	\includegraphics[scale=0.5]{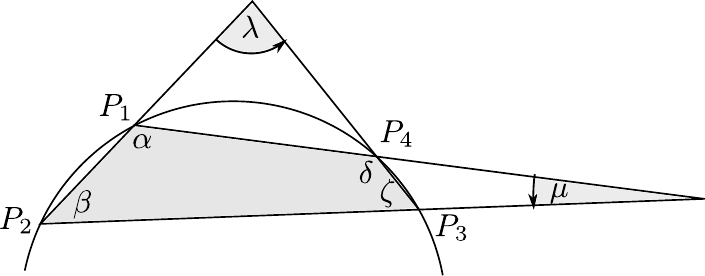}
	\caption{$Q$ and its angles.}
    \label{figLambdaMu}
  \end{minipage}
\quad
  \begin{minipage}[b]{0.45\textwidth}
    \centering
	\includegraphics[scale=0.3]{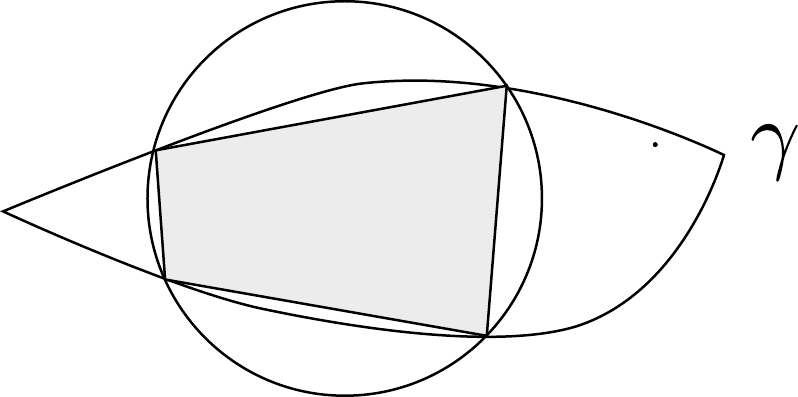}
	\caption{Example for Theorem~\ref{thmMainExtension}.}
    \label{figExampleThmExtension}
  \end{minipage}
\end{figure}

Consider a circular quadrilateral~$Q$.
We may and do assume that it is convex and positively oriented by relabeling its vertices in positively oriented fashion $P_1$, $P_2$, $P_3$, $P_4$.
Both pairs of opposite edges of $Q$ determine signed angles $\lambda = \measuredangle(\overrightarrow{P_1P_2},\overrightarrow{P_4P_3})$ and $\mu = \measuredangle(\overrightarrow{P_4P_1},\overrightarrow{P_3P_2})$, with the convention that $\lambda,\mu\in (-\pi,\pi)$, compare with Figure~\ref{figLambdaMu}.
Note that $\lambda$ and $\mu$ are zero, respectively, if and only if the corresponding pairs of opposite edges are parallel.
And if $\alpha$, $\beta$, $\zeta$, and $\delta$ denote the inner angles of $Q$, then
\[
\lambda = \alpha+\delta-\pi=\pi-\beta-\zeta \quad\textnormal{ and }\quad \mu = \zeta+\delta-\pi=\pi-\alpha-\beta.
\]
The following is a natural common extension of Theorems~\ref{thmMainTrapezoids} and~\ref{thmMainCircularQuads} for certain continuous convex Jordan curves and circular quadrilaterals, see Figure~\ref{figExampleThmExtension}.

\begin{theorem}[Common generalization]
\label{thmMainExtension}
Let $Q$ be a circular quadrilateral with signed angles $\lambda$ and $\mu$ as above.
Suppose $\gamma$ is a (continuous) convex Jordan curve all whose inner angles have size larger than $\min(|\lambda|,|\mu|)$.
Then $\gamma$ inscribes~$Q$.
\end{theorem}

The condition on the inner angles of $\gamma$ is only non-trivial at singular points of~$\gamma$, since $|\lambda|,|\mu|<\pi$.
In particular the angle condition is empty if $\gamma$ is $C^1$, and Theorem~\ref{thmMainCircularQuads} follows as a corollary.

As another special case, notice that if $Q$ is an isosceles trapezoid, then $\min(|\lambda|,|\mu|)=0$, which makes the angle condition again trivially fulfilled, and Theorem~\ref{thmMainTrapezoids} follows as a second corollary.

\paragraph{Related questions.}
There is a beautiful zoo of related theorems and open problems.
For example the reader may wonder about inscribed triangles in continuous curves (there are many, see Nielsen~\cite{Nie92trianglesInscribedInCurves}), or about inscribed pentagons (generically not possible, as the degree of freedom is one less then the number of equations).
We refer to various accounts on the history of inscribing and circumscribing problems, see Klee and Wagon~\cite[Problem 11]{KlWa96problemsInPlaneGeomAndNumberTh}, Nielsen~\cite{Nie10webpage}, Denne \cite{Den07squarePeg}, Karasev~\cite[2.6, 4.6]{Kar08topMeth}, Pak~\cite[I.3, I.4]{Pak09discreteAndPolyhedralGeometry}, M.~\cite{Mat12surveyOnSquarePeg}.

\paragraph{Basic ideas and outline.}
In the smooth case, we follow Karasev~\cite{KaVo10MakeevsConj}.
Given $Q=P_1P_2P_3P_4$, first one considers the set of inscribed triangles similar to $P_1P_2P_3$.
For generic curves $\gamma$, this set forms a one-dimensional manifold that winds around $\gamma$ exactly once, i.e.\ each of the three vertices circumscribe the interior of $\gamma$ once, see Section~\ref{secInscribedTriangles}.
Karasev's area argument then yields that the traced fourth vertex will circumscribe a region with the same signed area (see Corollary~\ref{corAreaOfP4}).
We will argue that if $\gamma$ does not inscribe $Q$, then this trace can be assumed to lie in the exterior of $\gamma$, going around $\gamma$ exactly once, and being injective, which yields a contradiction to the area argument.
The major new step here is to prove the injectivity of the trace, which is done in Section~\ref{secInjectivityOfP4}.

In the continuous case, two new problems arise:
Genericity of $\gamma$ and the corresponding approximation argument are harder to establish, which is a technical problem.
Furthermore, there is a new conceptual difficulty, namely that the inscribed triangles may become degenerate in a natural way, and at these singular points the traced fourth vertex may swap the sides of $\gamma$ without giving rise to a proper inscribed quadrilateral.
In many similar situations one would need to give up or find another approach (e.g.\ Toeplitz' inscribed square problem).
In our setting it turns out that after a more detailed analysis of these degenerate side changes in Section~\ref{secSingularPoints} we can actually use them to our advantage.
With inscribing problems it is often the case that the more complicated curves become, the more objects are inscribed, but to prove the existence of just a single one of them becomes harder (e.g.\ in the above two questions).
In our setting it seems to be quite the opposite.
We can even find a lower bound for the number of inscribed~$Q$'s, which can be tight even if the number of inscribed $Q$'s is large, see Theorem~\ref{thmQuantitativeThmMainExtension}.

Furthermore, the latter analysis can be used to show that inscribing $\Qtrap$ into $\CInf$ is equally difficult as inscribing them into the class $\COnePW$ of piecewise $C^1$ Jordan curves, see Section~\ref{secTrapezoidsOnNonConvexCurves}.

\paragraph{Notation.}

We say that two polygons $P_1P_2\ldots P_n$ and $Q_1Q_2\ldots Q_n$ in $\RR^2$ are \emph{similar} to each other if there is an orientation-preserving similarity transformation $\sigma$ (a composition of translations, rotations and scalings) such $Q_i=\sigma(P_i)$ ($i=1, \ldots, n$).

Throughout the paper, `smooth' means $C^\infty$.
We usually identify a parametrized curve $\gamma:S^1\to\RR^2$ with its image $\gamma(S^1)$ in order to simplify terminology.
We may and do assume that $\gamma$ goes in the positive sense around its interior.
Saying that $\gamma$ is $C^1$ or $C^\infty$ for us also includes that $\gamma$ needs to be regular.

We call a convex polygon $P_1P_2\ldots P_n$ positively oriented if $P_1$, \ldots, $P_n$ lie counter-clockwise around the boundary of the polygon.

Circular quadrilaterals may be self-intersecting (or ``skew''), in which case we can simply relabel the vertices in counter-clockwise order around the boundary of their convex hull, which makes the quadrilateral convex.
Inscribing either of them are equivalent tasks.
That is, it is enough to deal with positively oriented (and thus convex) circular quadrilaterals only.

%
%
%
%

\section{Inscribing the first three points}
\label{secInscribedTriangles}

Let us start with the easier smooth case.
Let $Q=P_1P_2P_3P_4$ be a circular quadrilateral with inner angles $\alpha,\beta,\gamma,\delta$.
For the sake of this paper we may assume that it is convex and positively oriented.
Furthermore we can cyclically permute the vertex labels to assure that $\zeta$ and $\delta$ are at least $\pi/2$, as $Q$ is circular.
This assumption will be crucial in Section~\ref{secInjectivityOfP4}.

\begin{figure}[htb]
  \centering
  \begin{minipage}[b]{0.45\textwidth}
	\centering
	\includegraphics[scale=0.4]{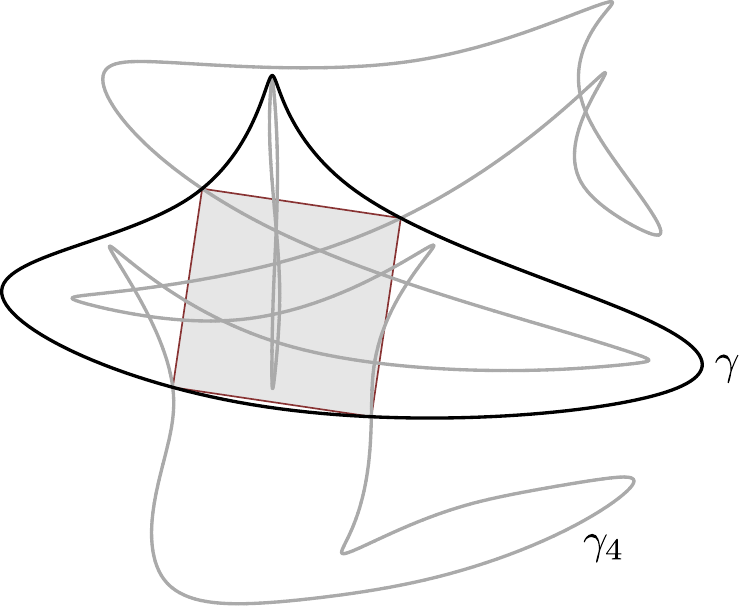}
	\caption{The trace $\gamma_4$ of the fourth vertex for a curve~$\gamma$ (exact drawing). Note that this exemplary curve is not convex.}
    \label{figTraceOGamma4}
  \end{minipage}
\quad
  \begin{minipage}[b]{0.45\textwidth}
    \centering
	\includegraphics[scale=0.4]{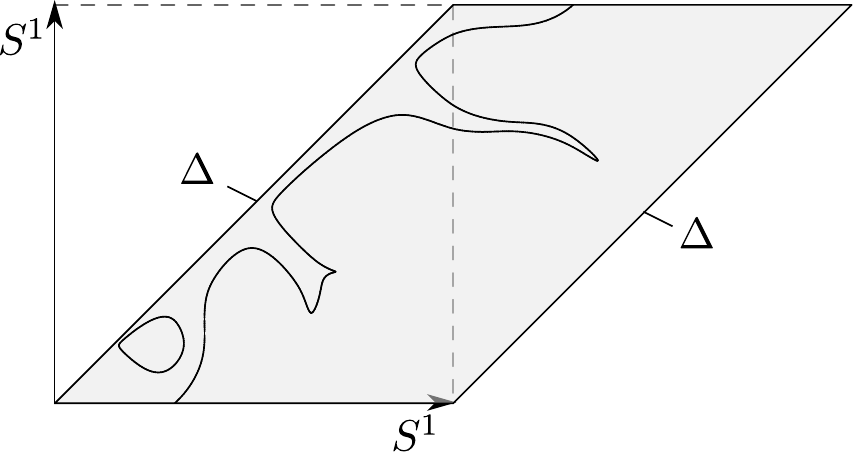}
	\caption{Projection of $Z_T\subset (S^1)^3$ to the first two coordinates for the curve from Figure~\ref{figTraceOGamma4} (exact drawing).}
    \label{figProjectionOfZT_exact}
  \end{minipage}
\end{figure}

Suppose $\gamma:S^1\to\RR^2$ is positively oriented $C^1$ convex Jordan curve.
We may deform $\gamma$ slightly (with respect to the $C^1$-metric) such that it becomes strictly convex and smooth.
If we can show that we can inscribe $Q$ into the deformed smooth strictly convex curve, the same follows for the original $\gamma$ by a limit argument, using $\gamma$ is $C^1$:
To any approximating smooth strictly convex curve we find an inscribed $Q$.
Making the approximation better in better (in the $C^1$-metric) yields a sequence of quadrilaterals, which by compactness has a convergent subsequence, whose limit cannot degenerate to a point as $\gamma$ is~$C^1$.

To any triangle $T = P_1P_2P_3$ we can consider the set $Z_T$ of all triangles $T'$ inscribed in $\gamma$ that are similar to~$T$, see Figure~\ref{figProjectionOfZT_exact} for the (non-convex) curve in Figure~\ref{figTraceOGamma4}.
This set has been studied topologically several times, also for more general polygons, see e.g.\ Meyerson~\cite{Mey80equilTriangles}, Wu~\cite{Wu04inscribingSmoothKnots}, Makeev~\cite{Mak05quadsInscribedInCurveAndVerticesOfCurve}, and \vrecica--\zivaljevic~\cite{VrZi08fultonMacPhersonCompCyclohedraAndPolygonalPegProblem}.

We consider $Z_T$ as a subset of the configuration space $(S^1)^3\wo\Delta$ ($\Delta$ always denotes a thin diagonal in this paper) which parametrizes all inscribed triangles in~$\gamma$.
As $Z_T$ can be defined via two equations, it can be written as a preimage $Z_T=f^{-1}(0)$ for some map $f:(S^1)^3\to\RR^2$.
Hence generically we expect that $Z_T$ is a one-dimensional proper submanifold.
The genericity can be achieved in various ways.
We choose to deform~$\gamma$ slightly with respect to the $C^1$-metric using `local bumps', keeping it strictly convex, where the amplitude of each bump depends on its own bounded real parameter.
Using the transversality theorem (see e.g.\ Guillemin and Pollack~\cite[p. 68]{GuPo10diffTop}) this makes $\gamma$ generic for any choice of amplitude vector outside a zero-set.
This method has the technical advantage that we did not deform the test-map, the curve itself becomes generic.
By an approximation argument as above, we thus may assume that $\gamma$ is not only smooth and strictly convex but also generic.

We claim that in case $\gamma$ is convex, $Z_T$ is topologically a circle; and even more is true:
For each angle $\alpha\in S^1=\ZZ/{2\pi}$ there is exactly one triangle parametrized by $Z_T$ whose first edge has angle $\alpha$ with the $x$-axis.
If there were more than one, these would be at least two inscribed triangles $T_1$ and $T_2$ that differ by a translation and a dilatation.
However then their six vertices cannot lie in strictly convex position. 
That for each $\alpha$ a corresponding inscribed triangle exists can be easily seen using an intermediate value theorem argument.
Or one computes directly the homology class that $Z_T$ represents, e.g.\ via a bordism argument deforming $\gamma$ to a simpler curve such as a circle.

Now for each such inscribed triangle $T' = P'_1P'_2P'_3$ we construct the fourth vertex $P'_4$ that makes $P'_1P'_2P'_3P'_4$ similar to the given~$Q$.
The trace of these points $P'_4$ is itself a closed curve $\gamma_4$ in the plane, although not necessarily simple, compare with Figure~\ref{figTraceOGamma4} for a non-convex curve~$\gamma$.
Now, each intersection point of $\gamma_4$ with $\gamma$ correspond to an inscribed quadrilateral $P'_1P'_2P'_3P'_4$ similar to the given~$Q$.
So assume that $\gamma_4 \cap \gamma = \emptyset$.
Then $\gamma_4$ stays inside of $\gamma$ or it stays outside.
We can restrict to the latter case by the following argument:
If we move a horizontal line parallelly from the bottom of $\gamma$ to the top, and at each time we call the intersection points $P_1''$ and $P_2''$, and if we construct corresponding points $P_3''$ and $P_4''$ to make $Q'' = P_1''P_2''P_3''P_4''$ similar to $Q$, then one of $P_3''$ or $P_4''$ will intersect $\gamma$ at last (if they do it simultaneously then we already are done with inscribing~$Q$).
If $P_3''$ comes last, then at that time, $T'' = P_1''P_2''P_3''$ lies in $Z_T$ and has $P_4''$ outside of~$\gamma$, which is the case we want to be in.
If $P_4''$ comes last, we simply relabel $1\leftrightarrow 2$, $3\leftrightarrow 4$, and reflect the plane and the orientation of $\gamma$ in order to arrive a positively oriented situation, and we arrived in the case, where the trace $\gamma_4$ of $P_4'$ stays outside of~$\gamma$.
The following lemma summarizes this.

\begin{lemma}
\label{lemSmoothCaseGenericSetting}
It is enough to prove Theorem~\ref{thmMainCircularQuads} for generic smooth strictly convex Jordan curves $\gamma$, and positively oriented circular quadrilaterals $Q$ with $\delta\geq\pi/2$ and whose trace $\gamma_4$ of $P_4'$ lies outside of~$\gamma$. 
\end{lemma}

For each $T'\in Z_T$, consider the intersection of the line segment between $P_2'$ and $P_4'$ with $\gamma$.
As $\gamma$ is convex, there is exactly one such intersection point $X$, except for $P_2'$ itself, and it moves continuously with $T'\in Z_T$.
We may consider this as a map $P_4'\mapsto X$.
Along $Z_T$, $P_2'$ winds once around $\gamma$ (possibly not in a monotone way), so does $X$, and thus the trace $\gamma_4$ of the fourth vertex $P_4'$ winds exactly once around $\gamma$ as well.

\section{On Karasev's and Tao's conserved integrals of motion}

Karasev~\cite{KaVo10MakeevsConj} proved that $\gamma_4$ circumscribes a region of signed area equal to the area of the interior of~$\gamma$.
Here, the signed area can be defined as one of the three equivalent integrals from Remark~\ref{remGreenIntegralsComparison}.
As a corollary he obtained the following theorem.
\begin{theorem}[Karasev]
Any smooth Jordan curve $\gamma$ either inscribes a given circular quadrilateral $Q=P_1P_2P_3P_4$, or it inscribes two copies of the triangle $P_1P_2P_3$ such that the two corresponding fourth vertices coincide.
\end{theorem}

His arguments behind this theorem are indeed the main ingredient for our proof of Theorem~\ref{thmMainCircularQuads}.

Tao~\cite{Tao17integrationApproachToToeplitzProblem} used a similar area argument in order to prove a new special case of Toeplitz' inscribed square problem, where the standard topological approach fails.

\begin{theorem}[Tao]
Let $f,g:[0,1]\to\RR^2$ be two $(1-\eps)$-Lipschitz functions whose graphs only intersect at $x=0,1$.
Then the curve formed by the two graphs inscribes a square.
\end{theorem}

And indeed his proof immediately generalizes to inscribed isosceles trapezoids, although one needs to assume a suitable smaller Lipschitz constant that depends on the angles of the trapezoid.

Whilst Karasev could use the fact that the four curves $\gamma_1$, $\gamma_2$, $\gamma_3$, $\gamma_4$ parametrizing the vertices of the quadrilateral in motion are closed, in Tao's situation they were not closed (at least in his application the path started and ended at quadrilaterals that were degenerate to a point). 

The two lemmas in this section simplify and extend Karasev's and Tao's conserved integrals of motion.
They work for arbitrary paths of circular quadrilaterals similar to the given one, which do not need to end where they started.
One hope is that the lemmas could be used in the future to help finding a proof of Makeev's conjecture that $\COne$ inscribes $\Qcirc$, for example by cutting the given curve into suitable pieces and applying the lemma to suitable $4$-tuples of these pieces.

\medskip

An affine dependence of points $P_1,\ldots,P_n$ in some $\RR$-vector space is a non-zero vector $(\lambda_1,\ldots,\lambda_n)\in\RR^n$ such that $\sum_i \lambda_i \begin{psmallmatrix} P_i \\ 1\end{psmallmatrix} = 0$, where $\begin{psmallmatrix} P_i \\ 1\end{psmallmatrix}$ denotes the projectivization of~$P_i$.
Any four points in the plane are affinely dependent.

\begin{lemma}[Area argument, complex version]
\label{lemAreaArgument}
Suppose $Q=P_1P_2P_3P_4$ is a circular quadrilateral.
Let $(\lambda_1,\ldots,\lambda_4)$ be an affine dependence of $P_1,\ldots,P_4$.
Let $\gamma_1,\ldots,\gamma_4: [t_0,t_1]\to\CC$ be four piecewise $C^1$-curves such that for each $t\in[t_0,t_1]$, the quadrilateral $\gamma_1(t)\gamma_2(t)\gamma_3(t)\gamma_4(t)$ is similar to~$Q$.
Then
\begin{equation}
\label{eqAreaArgument}
\sum_{i=1}^4 \lambda_i\int_{t_0}^{t_1} \gamma_i(t)\, d\overline{\gamma_i(t)} = 0.
\end{equation}
\end{lemma}

\begin{proof}[Proof of Lemma~\ref{lemAreaArgument}]
We proceed as in Karasev~\cite{KaVo10MakeevsConj}.
Let $O$ be the midpoint of $Q$ and $p_i = P_i - O\in\CC^\times$.
Let $\rho(t)\in\CC^\times$ denote the rotation-dilatation that sends $Q$ to a translate of $\gamma_1(t)\gamma_2(t)\gamma_3(t)\gamma_4(t)$, and let $o(T)$ denote the midpoint of $\gamma_1(t)\gamma_2(t)\gamma_3(t)\gamma_4(t)$.
Then clearly, $\gamma_i(t) = o(t) + \rho(t)p_i$.
Thus,
\begin{equation}
\label{eqgammadgammabar_expanded}
\gamma_i(t)d\overline{\gamma_i(t)} = o(t)d\overline{o(t)} +  p_i\rho(t)d\overline{o(t)} + \overline{p_i} o(t)d\overline{\rho(t)} + |p_i|^2 \rho(t)d\overline{\rho(t)}.
\end{equation}
If $r$ denotes the circumradius of $Q$, then
\[
\sum_i \lambda_i = 0,\quad
\sum_i \lambda_i p_i = 0,\quad
\sum_i \lambda_i \overline{p_i} = \overline{0} = 0,\quad
\sum_i \lambda_i |p_i|^2 = r^2\cdot 0 = 0.
\]
Thus, summing~\eqref{eqgammadgammabar_expanded} over $i=1, \ldots, 4$ with coefficients $\lambda_i$ yields
\begin{equation}
\label{eqAreaArgument_inTermsOfDifferentialForm}
\sum_{i=1}^4 \lambda_i\gamma_i(t)\,d\overline{\gamma_i(t)} = 0.
\end{equation}
Integrating this $1$-form over $t\in [t_0,t_1]$ yields~\eqref{eqAreaArgument}.
\end{proof}

In light of~\eqref{eqAreaArgument_inTermsOfDifferentialForm}, this seems to be in some sense the most natural formulation of the area argument.
The simplicity of the proof underlines that.
One possible caveat is that this talks about complex $1$-forms, so let us also discuss a version for real forms.

The $1$-forms $zd\bar z$ and $ydx$ on $\CC=\RR^2$ (with $z=x+iy$) are up to the factor $2i$ cohomologous (see below).
Therefore, Lemma~\ref{eqAreaArgument} can be rewritten in terms of $ydx$ as follows.

\begin{lemma}[Area argument, real version]
\label{lemAreaArgumentVia_ydx}
In the setting of Lemma~\ref{lemAreaArgument}, let $\rho_t = (\gamma_2(t)-\gamma_1(t))/(P_2-P_1)\in\CC^\times$ be the rotation-dilatation that sends $Q$ to a translated copy of $\gamma_1(t)\gamma_2(t)\gamma_3(t)\gamma_4(t)$.
Let $q$ be the quadratic form with matrix representation $\frac{1}{4}\begin{psmallmatrix}1&1\\1&-1\end{psmallmatrix}^t\left(\sum_{i=1}^4 \lambda_i P_iP_i^t\right)\begin{psmallmatrix}1&1\\1&-1\end{psmallmatrix}$.
Then
\begin{equation}
\label{eqAreaArgumentVia_ydx}
\sum_{i=1}^4\lambda_i\int_{\gamma_i} y\,dx = q(\rho_{t_1}) - q(\rho_{t_0}).
\end{equation}
\end{lemma}

\begin{remark}
1.) The `potential' $q$ in~\eqref{eqAreaArgumentVia_ydx} is a quadratic form on~$\RR^2$ of signature $(+,-)$.
Its two eigenvalues have opposite sign as the trace of $q$ is zero:
To show this, we may translate $Q$ to have its center at the origin, which keeps $q$ invariant.
Let $r$ be the radius of $Q$'s circumcircle.
Then $\Tr \sum_i\lambda_i P_iP_i^t = \sum_i\lambda_i|P_i|^2 = r^2\sum_i\lambda_i = 0$.
As $\begin{psmallmatrix}1&1\\1&-1\end{psmallmatrix}$ is $\sqrt{2}$ times an orthogonal matrix, the trace of $q$ is zero as well.

2.) Furthermore, the eigenvectors of $q$ are exactly the directions $v_\lambda,v_\mu$ of the angular bisectors of $\lambda$ and $\mu$.
Perhaps this has an elementary proof, but the author chose the brute-force algebraic way:
First one may assume that $v_\lambda,v_\mu$ are the standard basis vectors.
Then the coordinates of $P_1,\ldots,P_4$ satisfy a system of polynomial equations, giving rise to an ideal~$I$.
The statement about the eigenvectors of~$q$ is equivalent to say that $v_\lambda,v_\mu$ are isotropic vectors with respect to the quadratic form $\sum_i\lambda_i P_iP_i^t$, which translates into a polynomial equation in the coordinates of $P_1,\ldots,P_4$.
This polynomial is shown to lie in the ideal~$I$ using a Gr\"obner basis of~$I$, which was computed using SageMath~\cite{sage2017}, which in turn uses Singular~\cite{singular2016} for that task.

3.) Up to a scalar factor, the previous two points 1.) and 2.) uniquely describe $q$ geometrically.
\end{remark}

\begin{proof}[Proof of Lemma~\ref{lemAreaArgumentVia_ydx}]
Writing $z=x+iy$, we can expand and then collect terms
\begin{equation}
\label{eqzdzbar_expansion}
z d\overline{z} = xdx + ydy + iydx -i xdy = 2iydx + \tfrac{1}{2}d(x^2-2ixy+y^2),
\end{equation}
which shows that $zd\overline{z}$ and $2iydx$ differ only by a coboundary.
Summing up over all $i$ with coefficients $\lambda_i$ and on using~\eqref{eqAreaArgument} we obtain
\[
\sum_{i=1}^4\lambda_i\int_{\gamma_i} y\,dx = \frac{i}{4}\sum_{i=1}^4 \lambda_i\big(|\gamma_i(t_1)|^2  - 2i \gamma_i(t_1)_x\gamma_i(t_1)_y - |\gamma_i(t_0)|^2 + 2i \gamma_i(t_0)_x\gamma_i(t_0)y\big).
\]
Using $\sum_i \lambda_i|\gamma_i(t)|^2 = 0$ we can manipulate the right hand side further,
\begin{equation}
\label{eqSumLambdaiIntYdx_tmpExpression}
\sum_{i=1}^4\lambda_i\int_{\gamma_i} y\,dx = \frac{1}{4}\sum_{i=1}^4 \lambda_i\big(\langle\one,\gamma_i(t_1)\rangle^2 - \langle\one,\gamma_i(t_0)\rangle^2\big),
\end{equation}
where $\one = \begin{psmallmatrix}1 \\ 1\end{psmallmatrix}$.
We substitute $\gamma_i(t) = o(t) + \rho(t)p_i$ and use the affine dependence to get rid of the $o(t)$ summands and obtain $\sum_i \lambda_i \langle\one,\gamma_i(t)\rangle^2 = \sum_i \lambda_i \langle\one,\rho(t)p_i\rangle^2 = \sum_i\lambda_i \one^t M p_i p_i^t M^t \one$, where $M=\begin{psmallmatrix}a & -b \\ b & a\end{psmallmatrix}$ is the matrix representing the rotation-dilatation given by multiplication by $\rho(t) = a + ib$.
In the latter we can replace $p_i$ by $P_i$, as again via the affine dependence we see that the sum does not change.
Finally we write $M^t\one = \begin{psmallmatrix}a+b \\ a-b\end{psmallmatrix} = \begin{psmallmatrix}1 & 1 \\ 1 & -1\end{psmallmatrix}\rho(t)$.
This turns~\eqref{eqSumLambdaiIntYdx_tmpExpression} into the claimed~\eqref{eqAreaArgumentVia_ydx}.
\end{proof}

\begin{remark}
\label{remGreenIntegralsComparison}
If $\gamma:[t_0,t_1]\to\CC$ parametrizes a closed curve, then the signed area of the region circumscribed by $\gamma$ (counted with multiplicity) is given by Green's integrals $A=\int_{\gamma}x\, dy = -\int_{\gamma}y\, dx$.
With this in mind, integrating~\eqref{eqzdzbar_expansion} proves $\int_{t_0}^{t_1}\gamma\,d\overline{\gamma} = -2iA$.
\end{remark}

\begin{remark}
Taking as integrand $\gamma d\overline{\gamma}$ instead of $ydx$ has the advantage that the right hand side of~\eqref{eqAreaArgument} is simply~$0$ instead of the non-vanishing right hand side of~\eqref{eqAreaArgumentVia_ydx}, coming from the potential~$q$.
On the other hand, $ydx$ may have the advantage to be easier accessible geometrically, as it is immediately connected to areas.
\end{remark}

\begin{corollary}
\label{corAreaOfP4}
In the setting of Lemma~\ref{lemAreaArgument}, suppose that the $\gamma_i$ are closed curves.
If $\gamma_1$, $\gamma_2$ and $\gamma_3$ circumscribe regions of signed area $A$, then so does $\gamma_4$.
\end{corollary}

\begin{proof}
Let $\lambda$ be an affine dependence of $P_1,\ldots,P_4$.
Since any three vertices of a circular quadrilateral are affinely independent, $\lambda$ has no zero component.

If we put $A_i = -\int_{\gamma_{i}} y\,dx$, then either of Lemma~\ref{lemAreaArgument} and Lemma~\ref{lemAreaArgumentVia_ydx} implies $\sum_i \lambda_iA_i = 0$ since the curves are closed.
As $A_1=A_2=A_3=A$, $\sum \lambda_i = 0$, and $\lambda_i\neq 0$ for all $i$, $A_4$ needs to be equal to $A$ as well.
\end{proof}

\section{Injectivity of the fourth vertex' trace}
\label{secInjectivityOfP4}

In this section we finish the proof of Theorem~\ref{thmMainCircularQuads}.
In light of Lemma~\ref{lemSmoothCaseGenericSetting} and Corollary~\ref{corAreaOfP4}, it remains prove the following proposition.
Its proof relies heavily on the quadrilateral being cyclic.

\begin{proposition}
\label{propInjectivityOfP4}
Let $\gamma$ be a strictly convex smooth Jordan curve. Let $P_1P_2P_3P_4$ and $Q_1Q_2Q_3Q_4$ be two similar convex circular quadrilaterals with $P_4=Q_4$, such that the triangles $P_1P_2P_3$ and $Q_1Q_2Q_3$ lie counter-clockwise on $\gamma$, and such that $P_4=Q_4$ lies outside of $\gamma$, and such that the inner angle at $P_4$ is at least~$\pi/2$.
Then $P_1P_2P_3P_4 = Q_1Q_2Q_3Q_4$.
\end{proposition}

It reminds of math competition type problems. Indeed, it could be reformulated without mentioning~$\gamma$ at all, just assuming that $P_1,P_2,P_3,Q_1,Q_2,Q_3$ are in convex position but $P_4=Q_4$ lies outside their convex hull.

Before proving this proposition we need a lemma about circular quadrilaterals.
Any two distinct points $A$, $B$ in the plane determine a directed line $\overrightarrow{AB}$.
We say that a point $X$ lies \textdef{to the right of} $\overrightarrow{AB}$ if it lies in the closed half-space bounded by the line $AB$ that lies on our right hand side when we look from $A$ to~$B$.


\begin{lemma}
Let $P_1P_2P_3P_4$ be a convex circular quadrilateral. For $1\leq i<j\leq 3$ let $\rho_{ij}$ denote the rotation-dilatation about $P_4$ that sends $P_i$ to $P_j$. Then
\begin{enumerate}[label=\alph*)]
\item $\rho_{13}(\overrightarrow{P_1P_2}) = \overrightarrow{P_2P_3}$.
\item $\rho_{12}(\overrightarrow{P_1P_3}) = \overrightarrow{P_2P_3}$. 
\item $\rho_{23}(\overrightarrow{P_1P_2}) = \overrightarrow{P_1P_3}$.
\end{enumerate}
\end{lemma}

\begin{proof}
Let $\{i,j,k\} = \{1,2,3\}$.
Then $\rho_{ij}(P_kP_i) = P_kP_j$ follows from combining $\rho(P_i) = P_j$ and $\measuredangle P_iP_4P_j = \measuredangle P_iP_kP_j \mod \pi$.
This proves the lemma up to the orientation issue.
Now, all lines in the lemma are oriented in such a way that they have $P_4$ on their left as $P_1P_2P_3P_4$ is positively oriented, and all $\rho_{ij}$ fix $P_4$ and preserve the orientation of the plane, therefore they also respect the orientations of the lines.
\end{proof}

\begin{figure}[htb]
  \centering
  \begin{minipage}[b]{0.45\textwidth}
	\centering
	\includegraphics[scale=0.7]{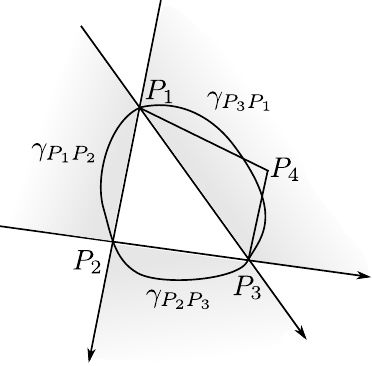}
	\caption{Possible regions for $\gamma_{P_iP_j}$.}
    \label{figGammaPiPj}
  \end{minipage}
\quad
  \begin{minipage}[b]{0.45\textwidth}
    \centering
	\includegraphics[scale=0.5]{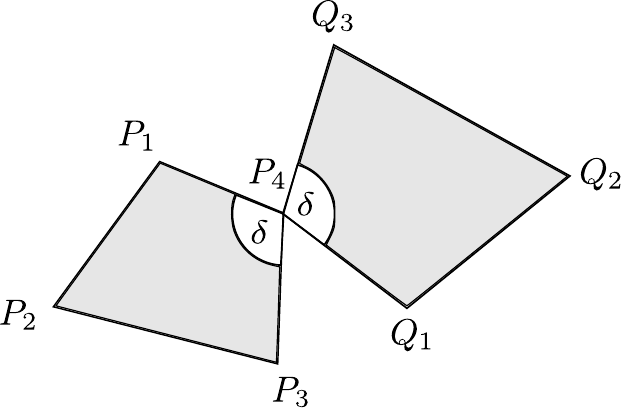}
	\caption{Case 3, $Q_1\in\gamma_{P_3P_1}$.}
    \label{figCase3}
  \end{minipage}
\end{figure}

\begin{lemma}
In the situation of Proposition~\ref{propInjectivityOfP4} the following equivalences hold.
\begin{enumerate}[label=\alph*)]
\item $Q_1$ lies to the right of $\overrightarrow{P_1P_2}$ if and only if $Q_3$ lies to the right of $\overrightarrow{P_2P_3}$.
\item $Q_1$ lies to the right of $\overrightarrow{P_1P_3}$ if and only if $Q_2$ lies to the right of $\overrightarrow{P_2P_3}$.
\item $Q_2$ lies to the right of $\overrightarrow{P_1P_2}$ if and only if $Q_3$ lies to the right of $\overrightarrow{P_1P_3}$.
\end{enumerate}
\end{lemma}

\begin{proof}
$Q_1Q_2Q_3Q_4$ is obtained from $P_1P_2P_3P_4$ via a rotation-dilatation about $P_4$, and any two rotation-dilatations about $P_4$ commute.
Thus $\rho_{ij}(Q_i) = Q_j$.
Therefore the lemma follows from the previous one using that rotation-dilatations preserve the orientation of~$\RR^2$.
\end{proof}

\begin{proof}[Proof of Proposition~\ref{propInjectivityOfP4}]
For two distinct points $A,B$ on $\gamma$, let $\gamma_{AB}$ denote the closed curve segment on $\gamma$ from $A$ to $B$ in counter-clockwise direction.
Then we get two decompositions
\[
\gamma = \gamma_{P_1P_2} \cup \gamma_{P_2P_3} \cup \gamma_{P_3P_1} = \gamma_{Q_1Q_2} \cup \gamma_{Q_2Q_3} \cup \gamma_{Q_3Q_1}.
\]

Note that each point $X\in\gamma_{P_iP_j}$ ($i\neq j$) lies to the right of~$\overrightarrow{P_iP_j}$, see Figure~\ref{figGammaPiPj}.

\paragraph{Case 1. $Q_1\in\gamma_{P_1P_2}$:}
Then $Q_1$ lies to the right of $\overrightarrow{P_1P_2}$, and thus by the lemma, $Q_3$ lies to the right of $\overrightarrow{P_2P_3}$, whence~$Q_3\in\gamma_{P_2P_3}$.

On the other hand, $Q_1$ lies to the right of $\overrightarrow{P_1P_3}$, and thus by the lemma, $Q_2$ lies to the right of~$\overrightarrow{P_2P_3}$.
Hence $Q_2\in\gamma_{P_2P_3}$ and thus $Q_2$ lies to the left of $\overrightarrow{P_1P_2}$, whence by the lemma, $Q_3$ lies to the left of~$\overrightarrow{P_1P_3}$.
Both restrictions on $Q_3$ only allow $Q_3=P_3$ and $P_1P_2P_3P_4 = Q_1Q_2Q_3Q_4$ follows.

From now on we may assume $Q_1\not\in \gamma_{P_1P_2}$, and by symmetry ($P\leftrightarrow Q$ and $1\leftrightarrow 3$),
\[
Q_1\not\in \gamma_{P_1P_2}, \quad P_1\not\in \gamma_{Q_1Q_2}, \quad
Q_3\not\in \gamma_{P_2P_3}, \quad P_3\not\in \gamma_{Q_2Q_3}.
\]

\paragraph{Case 2. $Q_1\in\gamma_{P_2P_3}$:}
Then $Q_1$ lies to the left of~$\overrightarrow{P_1P_3}$, hence by the lemma, $Q_2$ lies to the left of~$\overrightarrow{P_2P_3}$, and thus $Q_2\in\gamma_{P_2P_3}$.
As $P_3\not\in\gamma_{Q_2Q_3}$, together with $Q_2$ also $Q_3$ needs to lie in $\gamma_{P_2P_3}$.
Therefore, $Q_3$ lies to the right of~$\overrightarrow{P_2P_3}$, and hence by the lemma, $Q_1$ lies to the right of~$\overrightarrow{P_1P_2}$, whence $Q_1\in\gamma_{P_1P_2}$, which was already treated in the previous case.


\paragraph{Case 3. $Q_1\in\gamma_{P_3P_1}$:}
We may assume $Q_1\neq P_1$, otherwise the claim of the proposition follows.
As this case is the only remaining one, by symmetry we may assume
\[
Q_1\in\gamma_{P_3P_1}, \quad P_1\in\gamma_{Q_3Q_1}, \quad
Q_3\in\gamma_{P_3P_1}, \quad P_3\in\gamma_{Q_3Q_1}.
\]
This means not only that $Q_1$ and $Q_3$ lie on $\gamma_{P_3P_1}$, but also their order is determined:
In counter-clockwise order we see on $\gamma_{P_3P_1}$ the points $P_3,Q_1,Q_3,P_1$. (We allow that some of the points may coincide.)
As furthermore both triangles $P_1P_2P_3$ and $Q_1Q_2Q_3$ lie counter-clockwise on $\gamma$, this determines the cyclic order, in which all six of these points lie on $\gamma$, namely: $P_1,P_2,P_3,Q_1,Q_2,Q_3$ (up to cyclic permutation, and possibly $Q_1=P_3$ and/or $P_1=Q_3$).
This means (from $P_4$'th point of view, see Figure~\ref{figCase3}) that the two cones spanned by the angles $\angle P_1P_4P_3$ and $\angle Q_1Q_4Q_3$ with common apex $P_4=Q_4$ may at most have some boundary in common.
As the size of both angles was assumed to be at least $\pi/2$, it follows that $P_4=Q_4$ lies in the convex hull of $\{P_1,P_3,Q_1,Q_3\}\subset \gamma$, and thus not in the exterior of the convex curve $\gamma$, a contradiction!
\end{proof}

\begin{remark}
Without the angle restriction $\measuredangle P_1P_4P_3 \geq \pi/2$ one can indeed easily find two similar convex circular quadrilaterals $P_1P_2P_3P_4$ and $Q_1Q_2Q_3Q_4$ with $P_4=Q_4$, such that $P_1P_2P_3Q_1Q_2Q_3$ is a convex hexagon not containing $P_4=Q_4$. 
\end{remark}

\section{Singular curves}
\label{secSingularPoints}

Let $\gamma$ be a convex Jordan curve.
At each point $P\in\gamma$ we consider the inner angle $0<\alpha_P\leq\pi$ defined by $\alpha_P = \sup\measuredangle APB$, the supremum ranging over all points $A,B\in\gamma\wo\{P\}$.
Due to convexity, $P$ is a regular point of~$\gamma$ if and only if $\alpha_P=\pi$, otherwise it is a singular point of~$\gamma$.
We say that $\alpha_P$ is \emph{attained}, if this supremum is attained, i.e.\ if in a neighborhood of $P$, $\gamma$ looks like two straight line segments meeting at an angle $\alpha_P$.

The complementary angle at $P$ is $\alpha_P^c = \pi-\alpha_P$.
As the total curvature of $\gamma$ is $2\pi$, the sum of the complementary angles at all the singular points of $\gamma$ is at most $2\pi$.
This implies that there are at most countably many singular points, and for any $\eps>0$ there are at most finitely many singular points with $\alpha_P^c\geq\eps$, or equivalently, with $\alpha_P\leq \pi - \eps$.

Let $Q$ be a given circular quadrilateral with signed angles $\lambda$ and $\mu$ between their opposite edge pairs, as above Theorem~\ref{thmMainExtension}.
Then there are only finitely many singular points $S$ of $\gamma$ with inner angle $\alpha_S\leq\max(|\lambda|,|\mu|)$, let us call these the \emph{crucial} singular points, and we call $\alpha_S$ a \emph{crucial} angle.

They are crucial indeed, as they make the usual approximation argument break for two reasons.
\begin{enumerate}
\item If we smoothen $\gamma$ at a crucial singular point $S$ then this will introduce a tiny inscribed $Q$ close to~$S$ (unless $\alpha_S = \max(|\lambda|,|\mu|)$, in which case a case distinction is needed).
In the limit, there will be a sequence of such inscribed $Q$'s that converges to the quadrilateral that is degenerate to~$SSSS$.
\item When we trace inscribed triangles $T'$ similar to $T=P_1P_2P_3$, then they may also run into a crucial singular point $S$ and come out again in a different fashion, which is a priori not a serious problem.
The problem is that the trace $\gamma_4$ of the fourth vertex can change sides of $\gamma$, namely exactly when $T'$ is degenerate to~$SSS$.
Here, the area argument that worked for smooth curves would break, as the degenerate quadrilateral at $S$ does not count as an inscribed $Q$.
\end{enumerate}

If there is no crucial singular point, we can indeed simply approximate $\gamma$ by a smooth convex curve, reducing the problem to Theorem~\ref{thmMainCircularQuads}, and the limit argument works.

Let $S_1,\ldots,S_n$, be the crucial singular points, and $\alpha_{S_i}$ their inner angles.

\subsection{Reduction to a generic setting}
\label{secReductionToGenericSetting_extension}
Let $Q=P_1P_2P_3P_4$ be a circular quadrilateral with $\delta\geq \pi/2$, and let $\gamma$ be a convex Jordan curve.
The inner angles of $T=P_1P_2P_3$ are denoted by $\alpha_2,\beta,\zeta_2$.
Reducing to a generic setting is cumbersome for the above mentioned reasons.
We have to make sure that the approximation keeps the essential features of the curve such that we can easily study the neighborhoods of crucial singular points, and that the limit argument works (i.e.\ that finding a solution for each approximation yields a non-degenerate solution for the given curve).
Depending on the taste of the reader, we offer two different ways both leading to a useful generic approximation of $\gamma$, either a piecewise smooth one, or a piecewise linear one.
The author usually prefers smooth settings, however here the discrete one might indeed be less technical.

\subsubsection{Piecewise smooth approximation}
We proceed in three steps.
\begin{enumerate}
\item In case some $\alpha_{S_i} \in \{|\lambda|,|\mu|\}$, there are two possibilities:
\begin{enumerate}
\item If $\alpha_{S_i}$ is attained (i.e.\ the supremum in the definition of $\alpha_{S_i}$ is attained) then a neighborhood of $S$ looks like two line segments meeting at an angle~$\alpha_{S_i}$.
In that neighborhood, infinitely many copies of $Q$'s are inscribed.
\item If $\alpha_{S_i}$ is not attained, then we can deform the curve locally around $S_i$ making the inner angle slightly smaller and such that this smaller angle is attained.
\end{enumerate}
Now, no inner angle is equal to $|\lambda|$ or $|\mu|$.

\item We deform the curve  
such that it stays convex and is smooth away from the crucial singular points, and such that the crucial inner angles $\alpha_{S_i}$ are attained and do not belong to belong to $\{\alpha_2,\beta,\zeta_2\}$:
\begin{enumerate}
\item
If $\alpha_{S_i}$ is crucial, we replace a neighborhood of $\gamma$ around $S_i$ by two line segments that meet at some crucial angle close to $\alpha_{S_i}$ and not in $\{\alpha_2,\beta,\zeta_2\}$.
\item
If $\alpha_{S_i}$ is non-crucial, we replace a neighborhood of $\gamma$ around $S_i$ by a smooth arc.
\item
The remainder of $\gamma$ is deformed slightly in the $C^1$-sense to make $\gamma$ convex and smooth away from the crucial angles.
\end{enumerate}
This makes us easily understand the set $Z_T$ of inscribed triangles similar to~$T$ in the vicinity of singular points~$S_i$ (i.e.\ those triangles with all three vertices in a small neighborhood of~$S_i$ that are similar to~$T$):
There $Z_T$ is a union of smooth paths, with gaps exactly where $T'$ becomes degenerate to $S_iS_iS_i$. One could extend~$Z_T$ at these points continuously.

\item The test-map $f:(S^{1})^3\wo\Delta\to\RR^2$ from Section~\ref{secInscribedTriangles} that measured $Z_T = f^{-1}(0)$ may not be transversal to~$0$.
To solve this, we could add local bumps to~$\gamma$ as in Section~\ref{secInscribedTriangles} (which is possible). 
Instead, let us simply deform $f$ directly, as follows.
Around $\Delta$, $f$ is already transversal to~$0$ by the previous step.
So we deform $f$ only away from a neighborhood around~$\Delta$ by a suitable $\eps$-homotopy.
This makes its preimage $Z_T$ into a $1$-manifold, which is topologically a circle punctured at possibly some of the points $S_iS_iS_i$.
\end{enumerate}

\paragraph{Generic setting.} To summarize, we are now in the situation, where $\gamma$ is a convex Jordan curve with at most finitely many singular points, all of which are crucial, all of whose angles $\alpha_{S_i}$ are attained and not among $\{\alpha_2,\beta,\zeta_2,|\lambda|,|\mu|\}$.
And with the deformed test-map, $Z_T=f^{-1}(0)$ is a proper $1$-dimensional sub-manifold of $(S^1)^3\wo\Delta$, which parametrizes inscribed triangles that are up to some small error similar to $T$, and this error vanishes for small triangles.


\subsubsection{Piecewise linear approximation}
In case some $\alpha_{S_i} \in \{|\lambda|,|\mu|\}$, there are two possibilities:
\begin{enumerate}
\item If the supremum $\alpha_{S_i}$ is attained then a neighborhood of $S$ looks like two line segments meeting at an angle~$\alpha_{S_i}$.
In that neighborhood, infinitely many $Q$'s are inscribed.
\item If the supremum $\alpha_{S_i}$ is not attained, in what follows we will make sure to approximate this angle only from below (which can be done in general precisely because $\alpha_{S_i}$ is not attained).
\end{enumerate}


We construct a piecewise linear curve $\gamma_{PL}$ approximating $\gamma$ (in the $C^0$ sense) with the following properties:
It will be convex and piecewise linear.
Each inner angle of $\gamma_{PL}$ that approximates a crucial inner angle of $\gamma$ needs to be crucial as well. 
All other angles of $\gamma_{PL}$ must be non-crucial.
No inner angle of $\gamma_{PL}$ is allowed to be in $\{\alpha_2,\beta,\zeta_2,|\lambda|,|\mu|\}$.
So far this is actually not difficult to do.

Additionally we want that the set $Z_T$ of inscribed triangles similar to $T$ is a piecewise smooth $1$-manifold (in a generic way).
Here we use an algebraic trick.
If we pick a $3$-tuple $(e_1,e_2,e_3)$ of edges in $\gamma_{PL}$ and consider the triangles $T'=P_1'P_2'P_3'\in Z_T$ that have their $i$'th vertex on $e_i$ ($i=1,2,3$), then we see that they form a polytope:
We let $P_1'$ and $P_2'$ move freely on the lines extending~$e_1$ and~$e_2$, and see that the condition that $P_3'$ lies on the line extending $e_3$ is a linear equation. Furthermore, the restriction that $P_i$ lies on $e_i$ yields two linear inequalities (for each $i=1,2,3$).
So all we need to ensure is that these linear equations and inequalities are generic.
This can be achieved by choosing the vertices of $\gamma_{PL}$ in such a way that all its real coordinates are algebraically independent over the extension field $\QQ(\cos \alpha_2,\cos \beta, \cos \zeta_2,\cos\lambda,\cos\mu)$.
This is the promised algebraic trick.
We included the cosines such that none of the inner angles of $\gamma_{PL}$ is in $\{\alpha_2,\beta,\zeta_2,|\lambda|,|\mu|\}$. 

Thus in what follows we may work with $\gamma_{PL}$ in place of $\gamma$.

\subsection{Inscribed triangles in the neighborhood of singular points}
\label{secInscribedTrianglesAtSingularPoints}

Assume we are in the generic setting from above.
To simplify notation, we write $P_i$ instead of $P_i'$.

\begin{figure}[tbh]
  \centering
  \begin{minipage}[b]{0.88\textwidth}
	\centering
	\includegraphics[scale=0.7]{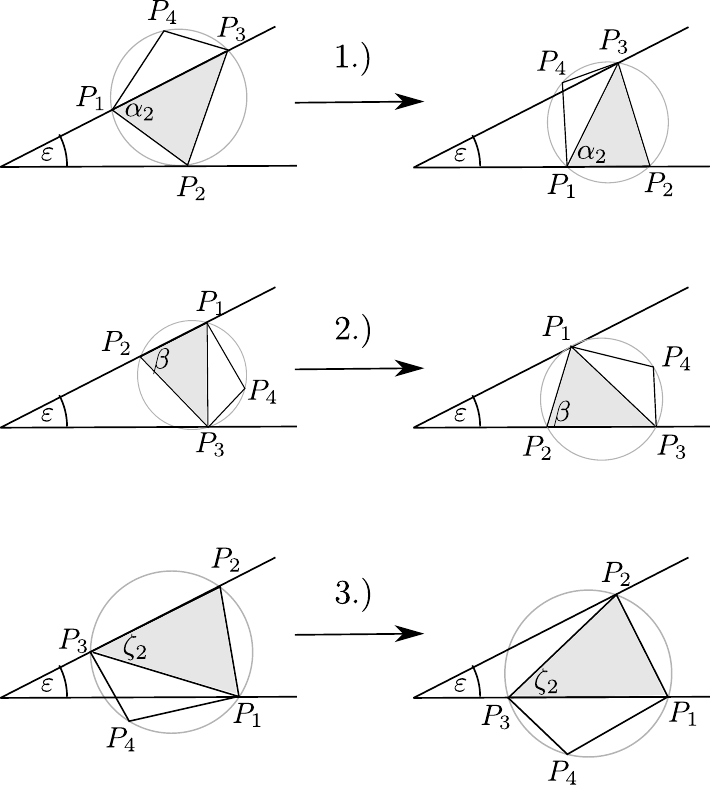}
	\caption{The three ways the curve $Z_T$ of inscribed triangles similar to $T$ can possibly move into a singular point, rotate the infinitesimal triangle by $\eps$, and come out again.}
    \label{figTrianglesAtSingularPoints}
  \end{minipage}
\end{figure}

Let $S$ be a singular point of $\gamma$ with inner angle $\alpha_S = \eps$ (crucial or not), and let $U$ be a sufficiently small neighborhood of $SSS\in (S^1)^3$.
$Z_T$ may visit $U$ up to three times, once for each inner angle of $T$ that is larger than~$\eps$.
To be precise, for each such inner angle, $Z_T\cap U$ has two components, namely one where the triangles run into the corner, and one where they come out of it.
In what follows, imagine that we connect these two ends with a one-dimensional family of imaginary infinitesimal triangles at $S$:
The component of $Z_T$ first runs into the corner $S$, becomes infinitesimally small, then it rotates on that spot counter-clockwise by the angle $\eps$, and finally it comes along $Z_T$ out of the corner again, see Figure~\ref{figTrianglesAtSingularPoints}.
All this can be made technically precise when working in the blowing-up of $(S^1)^3$ along~$\Delta$ (or in the Fulton--MacPherson compactification of $(S^1)^3\wo\Delta$).

Now let us see what happens with $P_4$ in each of the three cases, compare with Figure~\ref{figTrianglesAtSingularPoints}.

\begin{enumerate}[label=\arabic*.)]
\item This case occurs if and only if $\alpha_2>\eps$.
In the beginning, $P_4$ lies outside of $\gamma$.
At the end, $P_4$ lies outside if and only if $\lambda<\eps$.

\item This case occurs if and only if $\beta>\eps$.
In the beginning, $P_4$ lies outside of $\gamma$ if and only if $\lambda < -\eps$.
At the end, $P_4$ lies outside if and only if $\mu<-\eps$.
If $\lambda < -\eps$ or $\mu < -\eps$ (which implies $\beta>\eps$), $P_4$ will change sides of $\gamma$.

\item This case occurs if and only if $\zeta_2>\eps$.
In the beginning, $P_4$ lies outside of $\gamma$ if and only if $\mu < \eps$.
At the end, $P_4$ lies outside.
\end{enumerate}

This yields simple criteria for when $P_4$ changes sides of $\gamma$ when $Z_T$ passes a singular point $S$ with inner angle $\alpha_S = \eps$.
\begin{enumerate}
\item $P_4$ will change sides of $\gamma$ during motion 1.) at $S$ if and only if $\lambda>\eps$, as then $\alpha_2>\eps$ is automatically fulfilled.

\item $P_4$ will change sides of $\gamma$ during motion 2.) at $S$ if and only if $\lambda<-\eps$ or $\mu<-\eps$, as in this case $\beta>\eps$ is automatically fulfilled and both $\lambda,\mu<-\eps$ cannot happen as $\lambda+\mu = 2\delta-\pi \geq 0$.

\item $P_4$ will change sides of $\gamma$ during motion 3.) at $S$ if and only if $\mu>\eps$, as then $\zeta_2>\eps$ is automatically fulfilled.
\end{enumerate}

With this analysis we are ready to make use of these degeneracies.
We start with the proof of the general main theorem.

\begin{proof}[Proof of Theorem~\ref{thmMainExtension}]
In Section~\ref{secReductionToGenericSetting_extension} we argued what kind of genericity we can assume about~$\gamma$; as we may approximate non-generic curves by generic ones and use a limit argument.
Further we could assume that there is at least one crucial singular points as otherwise the proof for smooth curves can be used.
We will only consider the case $|\lambda|\geq|\mu|$; since the case $|\lambda|\leq|\mu|$ works analogously as the above criteria are essentially symmetric in $\lambda\leftrightarrow\mu$.

We will argue now how $P_4$ changes sides of $\gamma$ when $Z_T$ passes (a degenerate triangle at) a singular point $S$.
By the above criteria, whenever $Z_T$ passes a non-crucial singular point $S$, $P_4$ stays outside or stays inside.

More can happen at a crucial singular point $S$.
Here, $\alpha_S < |\lambda|$.
The assertion of the theorem implies $|\mu| < \alpha_S$.
As $\lambda+\mu = 2\delta-\pi\geq 0$, this can happen only if $\lambda>0$.
Thus, $|\mu| < \alpha_S < \lambda$.
Via the criteria above, we see that during the potential motions of type 2.) or 3.) at $S$, $P_4$ stays outside.
However there is a motion of type 1.) at $S$, and during that motion $P_4$ changes from the outside of $\gamma$ to the inside. 

These are all possibilities in which $P_4$ can change sides of $\gamma$ via a \emph{degenerate} inscribed~$Q$.
Thus, if $n$ is the number of crucial singular points $S$ of $\gamma$, then $P_4$ needs to go at least $n$ times back from the inside of $\gamma$ to the outside, and each time it yields a non-degenerate inscribed~$Q$.
\end{proof}

In fact, the proof can be easily made quantitative, which yields the following extension.

\begin{theorem}[Quantitative extension of Theorem~\ref{thmMainExtension}]
\label{thmQuantitativeThmMainExtension}
Let $Q$ be a circular quadrilateral with signed angles $\lambda$ and $\mu$ as above. Suppose $\gamma$ is a (continuous) convex Jordan curve all whose inner angles have size larger than $\min(|\lambda|,|\mu|)$. 
Let $n$ be the number of crucial singular points~$S$, i.e.\ those with inner angle $a_S\leq\max(|\lambda|,|\mu|)$.
Then $\gamma$ inscribes at least $\max(n,1)$ different copies of~$Q$.
\end{theorem}

\begin{proof}
We follow the proof of Theorem~\ref{thmMainExtension}.
Between any successive two of the $n$ times that $\gamma_4$ crosses $\gamma$ via a degenerate inscribed~$Q$, $\gamma_4$ needs to go back outside producing a non-degenerate inscribed~$Q$.
These events are separated from each other, use for example that using the discrete approximations the angle of $\overrightarrow{P_1P_2}$ with the $x$-axis is increasing (by how far depends only on the geometry of the original curve as well as on the $C^1$-distance of the original curve to its piecewise linear approximation).
Thus, in the limit we obtain $n$ different non-degenerate inscribed~$Q$'s for the given curve~$\gamma$.
\end{proof}

Given $Q$, the number $n$ of crucial singular points can bounded from above using the inequality $n(\pi-\max(|\lambda|,|\mu|))\leq \sum_{S\textnormal{ crucial}} \alpha_S^c\leq 2\pi$, which seems to bound the strength of Theorem~\ref{thmQuantitativeThmMainExtension}.

On the other hand, Theorem~\ref{thmQuantitativeThmMainExtension} can be tight for arbitrary large $n$:
Isosceles trapezoids have (after possibly relabeling the vertices) angles $\mu = 0$ and $0\leq \lambda < \pi$, and all such values for $\lambda$ are possible.
Now, a regular $n$-gon has $n$ inner angles of size $\alpha_{(n)} = \pi - 2\pi/n$.
For an isosceles trapezoid $Q$ with $\lambda > \alpha_{(n)}$, there are exactly $n$ ways to inscribe it in the regular $n$-gon, which matches the lower bound given in Theorem~\ref{thmQuantitativeThmMainExtension}.

\subsection{Inscribing isosceles trapezoids in non-convex curves}
\label{secTrapezoidsOnNonConvexCurves}
Akopyan asked (private communication) whether the implication Theorem~\ref{thmMainCircularQuads} $\impl$ Theorem~\ref{thmMainTrapezoids} proved in Section~\ref{secSingularPoints} works in the non-convex case as well.

To find a positive answer, let us restrict to the class $\COnePW$ of piecewise $C^1$ Jordan curves without cusps.
By this we mean curves $\gamma:S^1\to\RR^2$ that are $C^1$-regular along finitely many closed intervals that cover~$S^1$, and such that at the singular points of $\gamma$ there are no cusps (i.e.\ no inner or outer angles of size~$0$).

\begin{theorem}
\label{thmReductionOfCOnePWtoCInf_forInscribingQtrap}
$\COnePW$ inscribes $\Qtrap$ if and only if $\CInf$ inscribes $\Qtrap$.
\end{theorem}

One would conjecture that $\COnePW$ can be replaced by $\CZero$, but this is completely out of reach, as this would contain Toeplitz' inscribed square problem as a special case.

\begin{figure}[htb]
  \centering
  \begin{minipage}[b]{0.45\textwidth}
	\centering
	\includegraphics[scale=0.7]{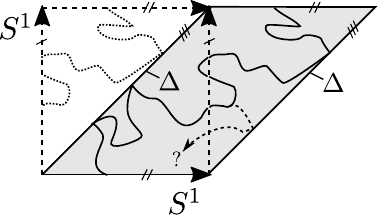}
	\caption{Projection of $Z_T^*$ to $(S^1)^2$. When bouncing off from $\Delta$, it cannot go `backwards'.}
    \label{figProjectionOfZT}
  \end{minipage}
\quad
  \begin{minipage}[b]{0.45\textwidth}
    \centering
	\includegraphics[scale=1.2]{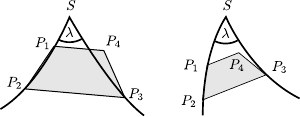}
	\caption{Singular point $S$ with inner or outer angle equal to $\lambda$ can become crucial or not.}
    \label{figInnerAngleLambda}
  \end{minipage}
\end{figure}

\begin{proof}[Proof of Theorem~\ref{thmReductionOfCOnePWtoCInf_forInscribingQtrap}]
We only need to prove the if-part.
Let $Q\in\Qtrap$, and we may assume $\lambda\geq 0$, $\mu=0$, $\delta\geq\pi/2$, and that $Q$ is positively oriented.
Let $\gamma\in\COnePW$. 
By assumption, any smooth approximation of $\gamma$ will inscribe of copy of~$Q$.

We extend the notion of the inner and outer angles in the obvious way to the non-convex setting: Let $\alpha_P\in(0,2\pi)$ denote the inner angle at $P\in\gamma$ (which shows towards the interior of $\gamma$), and $\alpha_P^\circ = 2\pi-\alpha_P$ the corresponding outer angle.
We call an (inner or outer) angle \emph{crucial} if it measures at most~$\lambda$.

If $\gamma$ has no crucial inner and outer angles, then we can use the standard limit argument to show that $\gamma$ inscribes a copy of~$Q$:
We approximate $\gamma$ suitably by a sequence of smooth Jordan curves $(\gamma_{(n)})_n$, such that the sizes of their inscribed copies of~$Q$ are uniformly bounded away from zero.
The latter is possible precisely because $\gamma$ has no crucial inner and outer angles.
By assumption, each $\gamma_{(n)}$ has an inscribed copy $Q_{(n)}$ of~$Q$, and by compactness, some subsequence converges to an inscribed copy of~$Q$ in~$\gamma$, which is non-degenerate due to the uniform lower bound on the sizes of~$Q_{(n)}$. 


We now come to the the case when $\gamma$ has a crucial inner or outer angle. 
This is the non-trivial case as smoothening a crucial angle will (usually) introduce a tiny inscribed copy of~$Q$, which would vanish in a limit argument into the corner.
First, we assume that $\gamma$ is as generic as needed with respect to the $C^1$-topology. 
%
%
%
As in the convex setting, we construct $Z_T$ and extend it with `infinitesimal' triangles at the singular points (as if the corners were infinitesimally smoothened) to make it into a closed $1$-manifold, sitting naturally in the blowing-up of $(S^1)^3\wo\Delta$ along~$\Delta$.
It still represents the same homology class in $H_1((S^1)^3;\ZZ)$, as is seen e.g.\ via a cobordism argument by deforming $\gamma$ into a strictly convex curve. 
However, $Z_T$ may have several connected components.
One can show that exactly one of these components represents the same homology class as~$\Delta$, i.e.\ it traces inscribed triangles each of whose vertices wind around $\gamma$ once, and all the other components are null-homologous; see Karasev~\cite{KaVo10MakeevsConj} for a formal proof.
Let $Z_T^0$ denote the former component of~$Z_T$, and $Z_T^{k}$ ($1\leq k\leq k_0$) the null-homologous ones.

As for convex curves the triangles~$T'=P_1P_2P_3$ traced by~$Z_T$ may run into singular points~$S\in\gamma$ and come out again, but now this is possible in two different ways, namely when one of the inner angles of $T$ is smaller than either the inner angle $\alpha_S$ (as before), or the outer angle $\alpha_S^o$ (the new case). 
Both cases are symmetric to each other, and in both settings we can use the analysis of Figure~\ref{figTrianglesAtSingularPoints}, except that they differ in up to two ways: 
\begin{enumerate}
\item The sides of the interior and the exterior of $\gamma$ are interchanged.
\item The direction of movement may be the opposite, i.e.\ for $\eps = \alpha_S^o$, the arrows in Figure~\ref{figTrianglesAtSingularPoints} may show the other way. 
This also depends on the orientation of~$Z_T$.
\end{enumerate} 
In fact one can show that if we give $Z_T$ one of the two possible preimage orientations, then that direction of movement will be opposite to Figure~\ref{figTrianglesAtSingularPoints}
exactly at outer angles; however we will only need a weaker statement (see the claim below).
As before, for each $T'$ we denote by $P_4 = P_4(T')$ the forth vertex that makes $P_1P_2P_3P_4$ similar to~$Q$.
For now we impose an additional genericity assumption on~$\gamma$ that none of its singular inner or outer angle is of size exactly~$\lambda$ (this will be justified in the last paragraph).
Then by the analysis of Section~\ref{secInscribedTrianglesAtSingularPoints}, during the passage of $Z_T$ through a singular vertex, $P_4$ will change sides with respect to $\gamma$ if and only if this is a motion of type $1$ at a crucial inner or outer angle.
Conversely, at each such angle, exactly one motion of type $1$ occurs (and possibly others of type $2$ and $3$).

Now consider one crucial inner or outer angle at some $S\in \gamma$, and let $Z_T^S$ be the component of $Z_T$ that passes $SSS$ in a motion of type~$1$.
\emph{We claim:} $Z_T^S$ can be oriented in such a way that each time $Z_T^S$ passes in that orientation some crucial (inner or outer) angle in a motion of type $1$, the fourth vertex $P_4$ moves from the outside to the inside of~$\gamma$.
As these are the only times where $Z_T^S$ passes a degenerate quadrilateral at which $P_4$ changes sides with respect to~$\gamma$, and since $P_4$ also has to move equally often from the inside of $\gamma$ to the outside, this claim proves the existence of an inscribed copy of~$Q$.


To prove the claim,
consider the projection $\pi_{12}:Z_T\to (S^1)^2$ given by the position of the first two vertices $P_1,P_2$ of the parametrized $T'$, see Figure~\ref{figProjectionOfZT}.
Let $Z_T^*:=\pi_{12}(Z_T^S)$, which is a closed path that may only touch the diagonal of $(S^1)^2$ without stepping over it.

\paragraph{Case 1: $Z_T^S = Z_T^0$.} Then $Z_T^*$ is homologous to the diagonal~$\Delta\subset(S^1)^2$, and we give $Z_T^*$ the orientation that corresponds to the standard orientation of~$\Delta$ (i.e.\ with tangent vectors $(1,1)$), and $Z_T^0$ the one that corresponds to it via~$\pi_{12}$.
Note that $Z_T^*$ does not self-intersect (except for possibly staying steady at some points $SS$ on the diagonal for some time) because $T'$ is determined by its edge~$P_1P_2$.
Therefore, whenever $Z_T^*$ touches~$\Delta$, it cannot `go back' (see the question mark in Figure~\ref{figProjectionOfZT}) as otherwise it would have to self-intersect by a Jordan curve theorem type argument.
We observe that at crucial inner angles the direction of motion when $Z_T^0$ passes $SSS$ is as in Figure~\ref{figTrianglesAtSingularPoints}, and at crucial outer angles it is the opposite.
This proves the claim in Case~1.

\paragraph{Case 2: $Z_T^S = Z_T^k$, $1\leq k\leq k_0$.} 
Note that $\pi_{12}(Z_T^0)$ cuts $(S^1)^2\wo\Delta$ into at least two connected components, and $Z_T^{*}$ must lies in one of them.
Thus, if $U\subset (S^1)^2$ is a small neighborhood of~$\Delta$, $Z_T^*$ intersects with only one of the two connected components of $U\wo\Delta$.
In other words, the only singular angles that $Z_T^S$ traverses are either all inner or all outer.
Moreover, by an analog ``no backwards'' argument as in Case 1, we see that the motions of passing a singular angle are always in the same direction as in Figure~\ref{figTrianglesAtSingularPoints} or always in the opposite direction (depending on the orientation of~$Z_T^S$).
The claim follows.

\smallskip

It still remains to discuss how $\gamma$ can be assumed to be generic.
This works as with convex curves (Section~\ref{secReductionToGenericSetting_extension}), but one crucial additional technical problem appears for singular points $S$ with $\alpha_S$ or $\alpha_S^o$ equal to~$\lambda$, see Figure~\ref{figInnerAngleLambda}.
If $\lambda=\alpha_S$ is an inner angle, we call the interior of $\gamma$ the \emph{$\lambda$-side} and the exterior the \emph{$\lambda^o$-side}; else $\lambda=\alpha_S^o$ is an outer angle and we swap these two notions.
In a suitably small neighborhood $U\subset S^1$ of $S$, the triangles $T'$ such that $P_2,P_1,S,P_3$ lie in this order on~$\gamma$ (or reversed) can be parametrized continuously: Near~$S$, for each $P_2$ there is exactly one such triangle, where for example $P_1$ can be obtained by intersecting $\gamma|_U$ with its own rotation about $P_2$ by the angle~$\beta$. (That this intersection point exists follows from $\alpha_2>\lambda$, and its uniqueness uses that $\gamma$ is composed of closed $C^1$-pieces and a mean value theorem type argument.)
Now consider the trace of $P_4$ when $T'$ approaches the degenerate triangle at~$S$.
If $P_4$ stays in the $\lambda$-side 
 of $\gamma$, we call $S$ crucial.
If $P_4$ stays in the $\lambda^o$-side of $\gamma$, we call $S$ non-crucial.
Otherwise $P_4$ intersects $\gamma$ on $T'$'s way towards $S$ and we are done.
Now, in the generic approximation of $\gamma$ that we construct, say $\gamma_{PL}$, we choose the inner/outer angle at $S$ to be strictly smaller or strictly larger than $\lambda$ depending on whether $S$ is crucial or not.
This keeps the trace of $P_4$ on the correct sides in the approximations, which avoids solutions that in the limit degenerate to~$S$.
\end{proof}

Two steps in the proof can be considered `lucky': 1.) We were able to use $2$-dimensional arguments of Jordan curve theorem type to show that $P_4$ can change sides of $\gamma$ only in one direction at points where~$Z_T$ degenerates.
The author is not aware of any other proofs in this area where technical difficulties for non-smooth curves arise in such a lopsided way that the theorem becomes trivially provable. 
2.) The inner or outer angles of size exactly $\lambda$, at which the notion of whether $P_4$ changes sides of $\gamma$ during motions of type $1$ may not be well-defined, can be deformed without negatively affecting the limit argument, as for the sake of $Z_T$ such angles are still generic.

\begin{remark}[Analogue of Theorem~\ref{thmMainExtension}]
\label{remReductionCOnePWtoCInf_forCertainQcirc}
Theorem~\ref{thmReductionOfCOnePWtoCInf_forInscribingQtrap} holds as well for circular quadrilaterals~$Q$ if we restrict to curves whose inner and outer angles are larger than $\min(|\lambda|,|\mu|)$; the proof is the same.
\end{remark}

Finally we can combine Theorem~\ref{thmReductionOfCOnePWtoCInf_forInscribingQtrap} with the recent result of Greene--Lobb~\cite{GreeneLobb20circularQuads} that~$\CInf$ inscribes~$\Qcirc$, and we obtain another positive partial answer for Question~\ref{quJ0inscribesQtrap}:

\begin{corollary}[{Assuming~\cite{GreeneLobb20circularQuads}}]
\label{corJ1pw_inscribes_trap}
$\COnePW$ inscribes~$\Qtrap$.
\end{corollary}

More generally, using Remark~\ref{remReductionCOnePWtoCInf_forCertainQcirc} in place of Theorem~\ref{thmReductionOfCOnePWtoCInf_forInscribingQtrap}, we obtain:

\begin{corollary}[{Assuming~\cite{GreeneLobb20circularQuads}}]
\label{corJ1pw_with_large_angles_inscribes_circ}
Any circular quadrilateral $Q$, with signed angles $\lambda$ and $\mu$ as above, can be inscribed into any $\gamma\in \COnePW$ whose inner and outer angles are all larger than $\min(|\lambda|,|\mu|)$.
\end{corollary}

The lower bound $\min(|\lambda|,|\mu|)$ is best possible, as can be seen by taking as $\gamma$ the union of two congruent circular arcs that meet at their endpoints in a given angle.


\paragraph{Acknowledgement.}
I wish to thank Arseniy Akopyan, Sergey Avvakumov, Roman Karasev and Sebastian Matschke for valuable correspondence.
In particular one of Akopyan's questions led to Section~\ref{secTrapezoidsOnNonConvexCurves}.
The plane geometry software Cinderella~\cite{RichterGebertKortenkamp12cinderella2manual} was a useful visualization tool when finding Proposition~\ref{propInjectivityOfP4} and its proof.
Moreover, I thank the anonymous referee for valuable comments.
This research was supported by the Initiative d'excellence de l'Universit\'e de Bordeaux (IdEx) and by Simons Foundation grant \#{}550023 at Boston University.

\small

\bibliographystyle{plain}

\bibliography{../../mybib07}

\begin{thebibliography}{10}

\bibitem{AkopyanAvvakumov17CyclicQuads}
Arseniy Akopyan and Sergey Avvakumov.
\newblock Any cyclic quadrilateral can be inscribed in any closed convex smooth
  curve.
\newblock {\href{http://arxiv.org/abs/1712.10205}{arXiv:1712.10205}}, 2017.

\bibitem{CDM11squarePeg}
Jason Cantarella, Elizabeth Denne, and John McCleary.
\newblock Transversality in {C}onfiguration {S}paces and the ``{S}quare-{P}eg''
  {P}roblem.
\newblock {\href{http://arxiv.org/abs/1402.6174}{arXiv:1402.6174}}, 2014.

\bibitem{Chr50kvadrat}
Carl~Marius Christensen.
\newblock A square inscribed in a convex figure (in {D}anish).
\newblock {\em Matematisk Tidsskrift B}, 1950:22--26, 1950.

\bibitem{singular2016}
Wolfram Decker, Gert-Martin Greuel, Gerhard Pfister, and Hans Sch\"onemann.
\newblock {Singular} {4-1-0} --- {A} computer algebra system for polynomial
  computations.
\newblock
  \href{https://www.singular.uni-kl.de}{https://www.singular.uni-kl.de}, 2016.

\bibitem{Den07squarePeg}
Elizabeth Denne.
\newblock Inscribed squares: Denne speaks.
\newblock
  {\href{http://quomodocumque.wordpress.com/2007/08/31/inscribed-squares-denne-speaks/}{http://quomodocumque.wordpress.com/2007/08/31/\-inscribed-squares-denne-speaks/}},
  2007.
\newblock Guest post on Jordan S. Ellenberg's blog {\emph{Quomodocumque}}.

\bibitem{sage2017}
The~Sage Developers.
\newblock {\em {S}age{M}ath ({V}ersion 8.0)}, 2017.
\newblock \href{http://www.sagemath.org}{http://www.sagemath.org}.

\bibitem{Emc13squarePeg1}
Arnold Emch.
\newblock Some properties of closed convex curves in a plane.
\newblock {\em Amer. J. Math}, 35:407--412, 1913.

\bibitem{Emc15squarePeg2}
Arnold Emch.
\newblock On the medians of a closed convex polygon.
\newblock {\em Amer. J. Math}, 37:19--28, 1915.

\bibitem{Emc16squarePeg3}
Arnold Emch.
\newblock On some properties of the medians of closed continuous curves formed
  by analytic arcs.
\newblock {\em Amer. J. Math.}, 38(1):6--18, 1916.

\bibitem{GreeneLobb20circularQuads}
Joshua~E. Greene and Andrew Lobb.
\newblock Cyclic quadrilaterals and smooth {J}ordan curves.
\newblock {\href{https://arxiv.org/abs/2011.05216}{arXiv:2011.05216}}, 2020.

\bibitem{GreeneLobb20rectangles}
Joshua~E. Greene and Andrew Lobb.
\newblock {The Rectangular Peg Problem}.
\newblock {\href{https://arxiv.org/abs/2005.09193}{arXiv:2005.09193}}, 2020.

\bibitem{GuPo10diffTop}
Victor Guillemin and Alan Pollack.
\newblock {\em Differential topology}.
\newblock Prentice Hall, 1974.

\bibitem{Hebbert14InscribedSquaresAndKinematicGeometry}
Clarence~M. Hebbert.
\newblock The inscribed and circumscribed squares of a quadrilateral and their
  significance in kinematic geometry.
\newblock {\em Ann. of Math. (2)}, 16(1-4):38--42, 1914/15.

\bibitem{vanHeijst14masterThesis}
Wouter {\noopsort{Heijst}{van Heijst}}, 2014.
\newblock Master thesis, in preparation.

\bibitem{Jer61inscribedSquares}
Richard~P. Jerrard.
\newblock Inscribed squares in plane curves.
\newblock {\em Trans. Amer. Math. Soc.}, 98:234--241, 1961.

\bibitem{Kar08topMeth}
Roman~N. Karasev.
\newblock Topological methods in combinatorial geometry.
\newblock {\em Russian Math. Surveys}, 63(6):1031--1078, 2008.

\bibitem{KaVo10MakeevsConj}
Roman~N. Karasev.
\newblock A note on {M}akeev's conjectures.
\newblock {\em J. Math. Sci.}, 212(5):521--526, 2016.

\bibitem{KlWa96problemsInPlaneGeomAndNumberTh}
Victor Klee and Stan Wagon.
\newblock {\em Old and new unsolved problems in plane geometry and number
  theory}.
\newblock Dolciani Mathematical Expositions. The Math. Ass. America, 1996.

\bibitem{Mak95quadsInscribedInClosedCurve}
Vladimir~V. Makeev.
\newblock On quadrangles inscribed in a closed curve.
\newblock {\em Math. Notes}, 57(1-2):91--93, 1995.

\bibitem{Mak05quadsInscribedInCurveAndVerticesOfCurve}
Vladimir~V. Makeev.
\newblock On quadrangles inscribed in a closed curve and the vertices of the
  curve.
\newblock {\em J. Math. Sci.}, 131(1):5395--5400, 2005.

\bibitem{Mat09squarePeg}
Benjamin Matschke.
\newblock On the {S}quare {P}eg {P}roblem and some relatives.
\newblock {\href{http://arxiv.org/abs/1001.0186}{arXiv:1001.0186}}, 2009.

\bibitem{Mat11phd}
Benjamin Matschke.
\newblock {\em Equivariant topology methods in discrete geometry}.
\newblock PhD thesis, Freie Universit{\"{a}}t Berlin, 2011.

\bibitem{Mat12surveyOnSquarePeg}
Benjamin Matschke.
\newblock A survey on the square peg problem.
\newblock {\em Notices Amer. Math. Soc.}, 61(4):346--352, 2014.

\bibitem{Mey80equilTriangles}
Mark~D. Meyerson.
\newblock Equilateral triangles and continuous curves.
\newblock {\em Polska Akademia Nauk. Fundamenta Mathematicae}, 110(1):1--9,
  1980.

\bibitem{Nie92trianglesInscribedInCurves}
Mark~J. Nielsen.
\newblock Triangles inscribed in simple closed curves.
\newblock {\em Geometriae Dedicata}, 43:291--297, 1992.

\bibitem{Nie10webpage}
Mark~J. Nielsen.
\newblock Web page on {F}igures {I}nscribed in {C}urves.
\newblock
  {\href{http://www.webpages.uidaho.edu/~markn/squares/}{http://www.webpages.uidaho.edu/$\sim$markn/\-squares/}},
  2000.

\bibitem{NielsenWright95rectanglesInscribedInSymmetricContinua}
Mark~J. Nielsen and Stephen~E. Wright.
\newblock Rectangles inscribed in symmetric continua.
\newblock {\em Geom. Dedicata}, 56(3):285--297, 1995.

\bibitem{Pak09discreteAndPolyhedralGeometry}
Igor Pak.
\newblock {L}ectures on {D}iscrete and {P}olyhedral {G}eometry.
\newblock
  {\href{http://www.math.ucla.edu/~pak/book.htm}{http://www.math.ucla.edu/$\sim$pak/book.htm}},
  2010.

\bibitem{OPT13noteOnToeplitzSquareProblem}
Ville~H. Pettersson, Helge~A. Tverberg, and Patric R.~J. {\"O}sterg{\aa}rd.
\newblock A note on {T}oeplitz' conjecture.
\newblock {\em Discrete Comput. Geom.}, 51(3):722--728, 2014.

\bibitem{RichterGebertKortenkamp12cinderella2manual}
J\"{u}rgen Richter-Gebert and Ulrich~H. Kortenkamp.
\newblock {\em {T}he {C}inderella.2 {M}anual}.
\newblock Springer-Verlag, 2012.

\bibitem{SaMa09inscribedSquaresDigitalPlane}
Feli{\'u} Sagols and Ra{\'u}l Mar{\'{\i}}n.
\newblock The inscribed square conjecture in the digital plane.
\newblock In {\em Combinatorial image analysis}, volume 5852 of {\em Lecture
  Notes in Comput. Sci.}, pages 411--424. Springer, 2009.

\bibitem{SaMa11inscribedSquares}
Feli{\'u} Sagols and Ra{\'u}l Mar{\'{\i}}n.
\newblock Two discrete versions of the inscribed square conjecture and some
  related problems.
\newblock {\em Theoret. Comput. Sci.}, 412(15):1301--1312, 2011.

\bibitem{Shn44geomPropClosedCurves}
Lev~G. Schnirelman.
\newblock On some geometric properties of closed curves.
\newblock {\em (in Russian) Usp. Mat. Nauk}, 10:34--44, 1944.
\newblock Available at
  {\href{http://ega-math.narod.ru/Nquant/Square.djv}{http://ega-math.narod.ru/Nquant/Square.djv}}.
  Posthumous reproduction and extension of the author's original article in
  \emph{Sbornik Rabot Matemati\v{c}eskogo Razdela Sekcii Estestvennyh i
  To\v{c}nyh Nauk Komakademii}, Moscow, 1929.

\bibitem{Str89inscribedSquares}
Walter~R. Stromquist.
\newblock Inscribed squares and square-like quadrilaterals in closed curves.
\newblock {\em Mathematika}, 36:187--197, 1989.

\bibitem{Tao17integrationApproachToToeplitzProblem}
Terence Tao.
\newblock An integration approach to the {T}oeplitz square peg problem.
\newblock {\em Forum Math. Sigma}, 5:e30, 63 pp, 2017.

\bibitem{Toe11aufgabenDerAnalysisSitus}
Otto Toeplitz.
\newblock {U}eber einige {A}ufgaben der {A}nalysis situs.
\newblock {\em Verhandlungen der Schweizerischen Naturforschenden Gesellschaft
  in Solothurn}, 4:197, 1911.

\bibitem{VrZi08fultonMacPhersonCompCyclohedraAndPolygonalPegProblem}
Sini{\v{s}}a Vre{\'{c}}ica and Rade~T. {\v{Z}}ivaljevi{\'{c}}.
\newblock {F}ulton--{M}ac{P}herson compactification, cyclohedra, and the
  polygonal pegs problem.
\newblock {\em Israel J. Math.}, 184(1):221--249, 2011.

\bibitem{Wu04inscribingSmoothKnots}
Ying-Qing Wu.
\newblock Inscribing smooth knots with regular polygons.
\newblock {\em Bull. London Math. Soc.}, 36(2):176--180, 2004.

\bibitem{Zin21konvexeGebilde}
Konrad Zindler.
\newblock {\"U}ber konvexe {G}ebilde.
\newblock {\em Monatshefte f\"ur Mathematik und Physik}, 31:25--56, 1921.

\end{thebibliography}

\end{document}